\documentclass[12pt]{article}
\usepackage[all]{xy}
\usepackage{amsfonts,amsmath,oldgerm,amssymb,amscd}
\newcommand{\N}{\mbox{$\mathcal{N}$}}
\newcommand{\M}{\mbox{$\mathcal{M}$}}
\newcommand{\G}{\mbox{$\mathcal{G}$}}
\newcommand{\Q}{\mbox{$\mathcal{Q}$}}

\newcommand{\J}{\mbox{$\mathcal{J}$}}
\newcommand{\h}{\mbox{$\mathcal{H}$}}

\newcommand{\F}{\mbox{$\mathcal{F}$}}
\newcommand{\I}{\mbox{$\mathcal{I}$}}

\newcommand{\ol}{\overline}
\newcommand{\R}{\mbox{$\mathbb R$}}
\newcommand{\Z}{\mbox{$\mathbb Z$}}

\newcommand{\q}{\mbox{$\mathbb Q$}}
\newcommand{\D}{\mbox{$\mathcal{D}$}}
\newcommand{\Hom}{\operatorname{Hom}}
\newcommand{\End}{\operatorname{End}}
\newcommand{\Aut}{\operatorname{Aut}}

\newtheorem{proof}{Proof}
\newtheorem{theorem}{Theorem}[section]
\newtheorem{proposition}[theorem]{Proposition}

\newtheorem{lemma}[theorem]{Lemma}
\newtheorem{remark}[theorem]{Remark}

\newtheorem{examples}[theorem]{Examples}

\newtheorem{notation}[theorem]{Notation}
\def\proof{\paragraph{Proof}}

\title{Anosov automorphisms on certain classes of nilmanifolds}
\author{Meera G. Mainkar\\
\small e-mail: meera@math.tifr.res.in}

\begin{document}

\maketitle

\begin{abstract}
We give a necessary and sufficient condition for  $k$-step
nilmanifolds  associated with graphs $(k \geq 3)$ to admit Anosov
automorphisms. We also prove nonexistence of Anosov automorphisms
on certain classes of $2$-step and $3$-step nilmanifolds.
\end{abstract}


\section{Introduction}\label{intro}
A well-known class of Anosov diffeomorphisms arises as follows.
 Let $N$ be a simply connected nilpotent Lie group and let $\Gamma$ be
a lattice in $N$; namely $\Gamma$ is a discrete subgroup  such
that $\Gamma \backslash N$ is compact. If $\tau$ is a hyperbolic
automorphism (see
 \S \ref{2-step} for the definition) of $N$
such that $\tau (\Gamma)=\Gamma$ then we get a diffeomorphism
$\overline \tau$ of $\Gamma \backslash N$, defined by $\overline
\tau (\Gamma x)= \Gamma \tau(x)$ for all $x\in N$, which is an
Anosov diffeomorphism of the compact nilmanifold $\Gamma
\backslash N$. Anosov diffeomorphisms arising in this way are
called {\it Anosov automorphisms of nilmanifolds}. Let  $K$ be a
finite group of automorphisms of $N$ and let  $\Gamma$ be  a
torsion free discrete cocompact subgroup of $K \ltimes N$.
  The $\Gamma$-action on $N$ is given by $(\tau , x).y = x \tau(y)$
where $\tau \in K$ and  $x,y \in N$.
 Now consider the quotient space $\Gamma \backslash N$ under the action of
$\Gamma$ on $N$. We call such a  compact manifold   $\Gamma
\backslash N$ an {\it infranilmanifold}. If $f$ is a hyperbolic
 automorphism  of $N$ such that $f$ normalises the subgroup $K$ in the group
 of automorphisms of $N$  and  $f(\Gamma) = \Gamma$  then
  $f$ induces a diffeomorphism  $\overline f$
 of the infranilmanifold  $\Gamma \backslash N$;  we call  such a $\overline f$ an
 {\it Anosov automorphism of an
  infranilmanifold $\Gamma \backslash N$}.

The only known examples of Anosov diffeomorphisms  are on  nilmanifolds and
 infranilmanifolds. It is conjectured that any Anosov diffeomorphism
is topologically conjugate to an Anosov automorphism of an
infranilmanifold.
 By a
result of A. Manning \cite{M} all Anosov diffeomorphisms on nilmanifolds are
topologically conjugate to Anosov automorphisms. This highlights the
question of
classifying   all
 compact nilmanifolds which admit Anosov automorphisms. Indeed it is
easy to see that not all of them do.  The first example (due to Borel) of a
 non-toral nilmanifold
admitting an Anosov automorphism was described by  S. Smale
\cite{S}. Later
 L. Auslander and J. Scheuneman
\cite{A-S} gave a  class of nilmanifolds admitting Anosov automorphisms.

 By a result of S. G. Dani (see  \cite{D1}), all
 nilmanifolds covered by free $k$-step nilpotent Lie groups on $n$
generators, with $k < n$, admit  Anosov automorphisms.
 There have been other recent
constructions of compact nilmanifolds with Anosov automorphisms (see
 \cite{D-M}, \cite{De-Des}, \cite{L}, \cite{L-W}
  and other references therein).

In this paper we associate a $k$-step nilmanifold $(k \geq 3)$
 with each graph, and
give a necessary and sufficient condition, in terms of the graph, for the
nilmanifold to admit Anosov automorphisms. We also prove some results on
nonexistence of Anosov automorphisms on certain $2$-step and $3$-step nilmanifolds.

\section{Preliminaries}\label{2-step}

In this section we recall some definitions and preliminaries
 concerning nilpotent
Lie groups and nilmanifolds. We also recall results concerning
 automorphisms of a
$2$-step nilmanifold associated with a graph (see \cite{D-M} for details).

 Let $N$
be a simply connected nilpotent Lie group and $\N$ be the
  Lie algebra of $N$, which is a nilpotent Lie algebra. Let $\Aut(N)$
 denote
the group of Lie  automorphisms of $N$. Let $\Aut(\N)$ denote the group of Lie
algebra automorphisms of $\N$. $\Aut(N)$ is isomorphic to the
 group $\Aut(\N)$, the
isomorphism is being given by $\tau \mapsto d\tau$, where $d\tau$ is the
differential of $\tau$.
 Let $\Gamma$ be a discrete subgroup of $N$ such that $\Gamma \backslash N$ admits
a finite $N$-invariant Borel measure. We call such a subgroup  a
{\it lattice} in $N$. As $N$ is a nilpotent Lie group, a discrete
subgroup $\Gamma$ is a lattice in $N$ if and only if $\Gamma
\backslash N$ is compact (see Theorem 2.1 in \cite{R}).

%\begin{definitions}\label{hyp}
 A {\it nilmanifold } is a quotient $\Gamma \backslash N$
 , where $N$ is a simply connected nilpotent Lie group and $\Gamma$ is
lattice in $N$.
 An automorphism $\sigma \in \Aut(\N)$ is said to be {\it hyperbolic } if
   all of its eigenvalues are of  modulus different from $1$.
 An automorphism $\tau \in \Aut(N)$ is said to be {\it hyperbolic } if
 all eigenvalues of the differential $d\tau$ are of
 modulus different from $1$.

%\end{definitions}

Now we recall the construction of the $2$-step nilmanifold
associated with a given graph and we recall some results about its
automorphism group (see \cite{D-M} for details).

Let $(S , E)$ be a finite simple graph, where $S$ is the set of
 vertices and $E$ is  the
set of edges. Let $V$ be a real vector space with $S$ as a basis.
Let $W$ be the subspace of  $\wedge^2V$
 spanned by
 $\{\alpha \wedge \beta :\alpha, \beta \in S, \alpha \beta \in
 E\}$, where  $\wedge^2V$ is the second exterior power of $V$. Let
 $\N= V \oplus W$.
We define the Lie bracket operation
 $[~,~]$ on $\N$
 as follows.   $[~,~]: \mathcal{N} \times
\mathcal{N} \rightarrow \mathcal{N}$
 is defined to be the unique bilinear map satisfying the following
 conditions:

i) for $\alpha, \beta \in S$, $[\alpha , \beta] = \alpha \wedge \beta$ if
$\alpha\beta \in E$ and $0$ otherwise;

ii) $[\alpha , \beta \wedge \gamma] = 0$ for all $\alpha, \beta, \gamma \in S$;

iii)  $[\alpha  \wedge \beta ,\gamma \wedge
  \delta] = 0$ for all $\alpha, \beta, \gamma, \delta \in S$.

We call $\N$ (defined as above) {\it the 2-step nilpotent Lie algebra
 associated to the graph $(S,E)$}.
%By using the Baker-Campbell-Hausdorff formula
%(see \S 2.15 in \cite{V}),
%$N$ can be realised as $\N$ with the multiplication
%defined by
%\[(v_1,w_1)\cdot (v_2,w_2)= (v_1+v_2, w_1+w_2 +\frac12 [v_1,v_2]),\]
%for all $v_1,v_2 \in V$ and $w_1,w_2 \in W$ (here $+$ denotes the addition in
%$\N$)
  Let $N$ be the simply connected Lie group with Lie algebra $\N$. Let $\Gamma$
be the subgroup of $N$ generated by exp$(S)$, where exp denotes
the
 exponential map.
   It can be seen that $\Gamma$ is a lattice in $N$.
A nilmanifold $\Gamma \backslash N$
   is called  {\it the $2$-step nilmanifold associated with the graph $(S,E)$}.

For any $\sigma \in S$ we define
$$\Omega' (\sigma)=\{\omega  \in S: \sigma \omega \in E\} \ \  \hbox{\rm
  and}\ \  \Omega (\sigma)=\Omega' (\sigma)\cup \{\sigma\}. $$
 Let $\sim$ be an equivalence relation on $S$ defined as follows:
 for $\alpha, \beta \in S$,
 $\alpha \sim \beta $ if either $\alpha =\beta$ or,  $\Omega'(\alpha) \subset
 \Omega(\beta)$ and $\Omega'(\beta) \subset \Omega(\alpha)$ (see
\cite{D-M}
for  details).
 Let $\{S_\lambda\}_{\lambda \in \Lambda}$ denote the set of all equivalence
 classes in $S$ with respect to the equivalence relation $\sim$, where $\Lambda$ is
 an index set. $S_\lambda$, $\lambda \in \Lambda$, are called the {\it
coherent components} of $S$.
 For each $\lambda \in \Lambda$, let $V_\lambda$ denote the subspace
 of  $V$ spanned by $S_\lambda$.

We recall some results (see \cite{D-M}):

\begin{theorem}\label{thm}
Let $(S,E)$ be a finite graph and let $\N = V \oplus W$ be the
 $2$-step  nilpotent Lie algebra associated with $(S,E)$ (notation
as above). Let $G$ denote the subgroup of $GL(V)$ consisting of
all restrictions, $\tau|V$, such that $\tau \in \Aut(\N)$ and
$\tau(V) = V$. Then $G$ is a Lie subgroup of $GL(V)$ and the
following conditions are satisfied:
%i) The Lie subalgebra $\G$ of  $G$ is given by $$\G=\D\oplus  \bigoplus_{(\alpha,
%\beta) \in \Theta(S,E)} \R E_{\alpha\beta}\,.$$

i) The connected component of the identity in $G$, which we denote be $G^0$,
   can be expressed as $(\prod_{\lambda \in \Lambda}\,GL^+(V_\lambda))\cdot M,$
where for each $\lambda \in \Lambda $, $GL^+(V_\lambda)$ denotes the subgroup of
$GL(V_\lambda)$ consisting of all the elements with positive determinant and $M$ is
a closed connected nilpotent normal subgroup of $G$.

ii)  The elements of $\Lambda$ can be arranged as $\lambda_1, \dots , \lambda_k$ so
that for all $j=1,\dots, k$,   $\bigoplus_{i\leq j}\,V_{\lambda_i}$ is invariant
under the
  action of $G^0$.

\end{theorem}

\begin{lemma}\label{lem:det}
Let $\lambda_1, \dots ,\lambda_k$ be an enumeration of $\Lambda$ such that
assertion~(ii) of Theorem~\ref{thm} holds. For each $j=1,\dots,k$ let
 $\N_j = (\oplus_{i\leq j}\,V_{\lambda_i})\oplus W$; also let
$\N_0=W$. Let $\tau$ be a Lie automorphism of $\N$ contained in
the connected component of the identity in $\Aut (\N)$. Then each
$\N_j$ is invariant under the action of $\tau$. Let $\Phi$ be the
(additive) subgroup of $\N$ generated by $S \cup \{\frac{1}{2}
(\alpha \wedge \beta : \alpha, \beta \in S, \alpha\beta \in E \}$.
 If $\tau (\Phi)=\Phi$ then for all $j=1,\dots , k$ the determinant
of the action of $\tau$ on $\N_j$ is $\pm 1$.
\end{lemma}

\section{$k$-step nilmanifold associated with the
graph}\label{k-stepdef}

In this section we associate a $k$-step $(k \geq 3)$ nilmanifold
 (i.e.~covered by a $k$-step simply connected nilpotent Lie group)
with every graph and  we give a necessary and sufficient condition for such
nilmanifolds to admit an Anosov automorphism.

 Starting with a graph $(S, E)$ we  define a $k$-step $(k \geq 3)$
 nilpotent Lie algebra  as follows.
 Let $(S,E)$ be a finite graph, where $S$ is the
 set of vertices and $E$ is the set of edges. Suppose $\N$ denotes
 the 2-step nilpotent Lie
algebra
 associated with $(S,E)$ (see \S \ref{2-step}) i.e.
 $\N= V \oplus W$ where
$V$ is a vector space with $S$ as a basis and $W$ is  the subspace of $\wedge^2V$
 spanned by $\{\alpha \wedge \beta : \alpha, \beta \in S, \alpha \beta \in E
 \}$.

Let $\N_k(V)$ be a free $k$-step nilpotent Lie algebra on $V$
 (see \cite{A-S} for
 the definition).  We denote by  $\h_k$
the $k$-step
nilpotent Lie algebra  $ \N_k(V)/\J$, where
 $\J$ denotes an ideal of $\N_k(V)$ generated by all  elements $[\alpha,
\beta]$ such that $\alpha \beta$ is not an edge.  Let $N_K$ be the
simply connected nilpotent Lie group with Lie algebra $\h_k$.
  Suppose $\Phi_k$ is the
(additive) subgroup of $\h_k$ generated by the elements of the type
$[\alpha,[\beta,\dots]]$, where $\alpha, \beta, \ldots \in S$.
 %Let $N_K$ be the
%simply connected nilpotent Lie group with Lie algebra $\h_k$.
 Then there exists a $\Z$-subalgebra $\Phi_k^0$  of $\h_k$ which is
 contained in $\Phi_k$ such that $\Gamma_k = $ exp$(\Phi_k^0)$
 is a subgroup of $N_k$, and $\tau(\Phi_k) = \Phi_k$ if and only if
  $\tau(\Phi_k^0) = \Phi_k^0$ for any automorphism $\tau$ of
  $\h_k$
(see  \S2 in \cite{A-S}).
   We note that $\Gamma_k$ is a lattice in $N_k$.  We call a
 nilmanifold $\Gamma_k \backslash N_k$  {\it a
 $k$-step nilmanifold associated with the graph $(S,E)$}.

\section{Anosov automorphisms of $\Gamma_k \backslash N_k , k \geq 3$}\label{auto}
We give a necessary and sufficient condition for the nilmanifold
 $\Gamma_k \backslash N_k$ to
 admit an Anosov automorphism.

\begin{notation}\label{M^n}
{\rm  Suppose $\h$ is a $k$-step nilpotent Lie algebra. For any subset $M$ of $\h$
we denote $[M, M]$ by $M^1$, $[M, [M, M]]$ by $M^2$, and inductively we denote
  $[M, M^{n-1}]$ by $M^n$ for all $n$, $3 \leq n \leq k-1$.
}
\end{notation}

 Let $(S, E)$ be a graph and $\Gamma_k \backslash N_k$ be a
 $k$-step nilmanifold ($k \geq 2$) associated
with $(S, E)$. We refer \S \ref{2-step} and \S \ref{k-stepdef} for the notation.
\begin{remark}\label{aut}
\rm{  We note that any automorphism  of $\N$ can be extended to an automorphism of
$\h_k$.   The automorphism group
 $\Aut(\h_k)$ is the semidirect product of
$\Aut(\N)$ and a connected group. This can be seen by observing that $\Aut(\h_k)$ is
a semidirect product of $\Aut(\h_k/\h_k^{k-1})$ and $\Hom(V, \h_k^{k-1})$, and
$\Aut(\h_k/\h_k^{2})$ is the same as $\Aut(\N)$.
% Since $\Aut(\N)$ has finitely
%many connected components,  $\Aut(\h_k)$ has finitely many connected
%components.
}
\end{remark}

\begin{theorem}\label{kstep}
$\Gamma_k \backslash N_k$ admits an Anosov  automorphism if and
only if
 the following holds:

i) For every $\lambda, |S_\lambda| \geq 2$; and

ii) If $|S_\lambda| = l$, with  $2 \leq l \leq k$, and $\alpha,
 \beta\in S_\lambda$ then
 $\alpha \beta$ is not an edge.
\end{theorem}

\proof Suppose that for each $\lambda \in \Lambda$ (i) and (ii) hold. We will prove
that there exists a hyperbolic automorphism
%(see \ref{intro} for the definition).
 $\tau \in \Aut(\h_k)$ such that $\tau(\Phi_k) = \Phi_k$.

 We have
 $V=\oplus_{\lambda \in \Lambda}\,V_\lambda$ (see \S \ref{2-step} for notation).
 For each $\lambda \in \Lambda$ let
$\Phi_\lambda$ be the subgroup of $V_\lambda$ generated by $S_\lambda$. There exists
$g_\lambda \in   GL(V_\lambda)$ such that $g_\lambda(\Phi_\lambda) =\Phi_\lambda $
if and only if the matrix representing $g_\lambda$ with respect to the basis
$S_\lambda$ belongs
 to $GL(d_\lambda, \Z)$, where   $d_\lambda =|S_\lambda|$.
For each $\lambda \in \Lambda$ there exists a matrix $A_\lambda \in   GL(d_\lambda,
\Z)$
 with the eigenvalues $c_1, c_2, \ldots c_{d_\lambda}$ such that
 $|c_{i_1} c_{i_2} \cdots c_{i_r}| \neq 1,$  for all $ r, 1 \leq r
\leq $  min $(k, d_\lambda -1)$,  and for all
 $ i_1, i_2, \ldots, i_r \in \{1, 2, \ldots, d_\lambda \}. $
The existence of such elements can be proved by using a result of S. G. Dani
  (see Corollary 4.7 in \cite{D}).
 Let $g_\lambda$ denote the transformation from $GL(V_\lambda)$
 whose matrix  with respect to the basis $S_\lambda$ is $A_\lambda$.
By the above observation $g_\lambda(\Phi_\lambda) =\Phi_\lambda$.
 We choose natural numbers $j_\lambda$, $\lambda \in \Lambda$,  such
that $|\prod_{\lambda \in \Omega}
 ( c_{\lambda i_1}c_{\lambda i_2} \cdots
c_{\lambda i_{n_\lambda}})^{j_\lambda}| \neq 1  $  for all subsets
$\Omega$
 of $\Lambda$  such that $|\Omega| \geq 2$ and $2 \leq \sum_{\lambda \in \Omega}
n_\lambda \leq k$, where  $c_{\lambda i_j}$'s  are eigenvalues of
$g_\lambda$. Let $g \in GL(V)$ be the element whose restriction to
$V_\lambda$ is $g_\lambda^{j_\lambda}$, for each $\lambda \in
\Lambda$.

  There exists $\overline \tau \in \Aut(\N)$ such that $g$ is the
 restriction of $\overline \tau$ to $V$ (see Theorem \ref{thm}).
 We know that $ \overline \tau$ constructed as above is a hyperbolic automorphism of
$\N$. This can be seen from  the proof of the Theorem $1.1$ in \cite{D-M} and  the
hypothesis of the theorem.
 Let $\tau$ be an automorphism of $\h_k$ obtained by extending
$\overline \tau$.  We note that $\tau(\Phi_k) = \Phi_k$ by construction.
 We will prove that $\tau$ is  hyperbolic as a linear
transformation. Suppose if possible $\tau$ has an eigenvalue, say $c$, of absolute
value $1$. Then $c$ must be an eigenvalue of the restriction of $\tau$ to   $V^n$,
$3 \leq n \leq k$ (see Notation \ref{M^n}), since $\tau$ is hyperbolic on $\N$.

Now using the fact that $\tau(V_\lambda) = V_\lambda$ for all $\lambda \in \Lambda$
and recalling  the construction of $g$, we see that
 there exists $\lambda \in
\Lambda$ such that $|S_\lambda| = n$ and $V_\lambda^{n}$ is nonzero (see
 Notation \ref{M^n}). But by the condition in the
hypothesis $\alpha \beta$ is not an edge, for all $\alpha,
 \beta \in  S_\lambda$.   Hence
$[\alpha, \beta] = 0$, for all $\alpha,
 \beta \in  S_\lambda$. This contradiction shows that $\tau$ is
hyperbolic. Hence $\Gamma_k \backslash N_k$ admits an Anosov
automorphism.

Conversely suppose that $\Gamma_k \backslash N_k$ admits an Anosov
automorphism. Hence there exists $\tau \in \Aut(\h_k)$ such  that
$\tau(\Phi_k)= \Phi_k$ and $\tau$ is a hyperbolic linear
transformation. Let $\overline \tau \in \Aut(\N)$ denote an
automorphism of $\N$ induced by $\tau$. We can assume that
$\overline \tau (\Phi) = \Phi$, where $\Phi$ is the
 subgroup of $\N$ (with respect to addition) generated by the subset
$S \cup \{\frac12 (\alpha\wedge \beta):\alpha, \beta \in S, \alpha \beta \in E\}$.
As $\overline \tau$ is a hyperbolic linear transformation, $|S_\lambda| \geq 2$
  for every $\lambda$, and if
$|S_\lambda| = 2$ then $\alpha \beta$ is not an edge for $\alpha \beta \in
S_\lambda$ (see Theorem 1.1 in \cite{D-M}). We may assume that $\overline \tau$ is
contained in the connected component of the identity in $\Aut(\N)$ (see  Remark
\ref{aut}).
    Let $G$ denote the subgroup of $GL(V)$ consisting of all restrictions, $\tau|V$,
  such that $\tau \in \Aut(\N)$ and  $\tau(V) = V$.   We write the elements of
$\Lambda$ as $\lambda_1, \lambda_2, \dots, \lambda_m$ such that for all $j= 1,\dots,
m$, $\bigoplus_{i\leq j}\,V_{\lambda_i}$ is invariant under the  action of $G^0$,
 where $G^0$ is the connected component of identity in $G$ (see Theorem  \ref{thm}).
   Now suppose there exists $\lambda \in \Lambda$ such that $|S_\lambda| = l$, $3
\leq l \leq k$  and $\alpha \beta  \in E$ for all $\alpha, \beta \in S_\lambda$. Let
$j, 1 \leq j \leq m$, be such that $\lambda = \lambda_j$.
 Consider the induced action of $\overline \tau$ on $\N_j/\N_{j-1}$,
where $\N_j =  (\oplus_{i\leq j}\,V_{\lambda_i})\oplus W$. We note that  each $\N_j$
is invariant under the action of $\overline \tau$
 (see  Lemma \ref{lem:det}). As
the determinant of the induced action  of $\overline \tau$ on $\N_j/\N_{j-1}$
 is $\pm 1$, the product of the eigenvalues
$\theta_1, \theta_2, \dots, \theta_l$ of the induced action is
$\pm 1$. Since the action is hyperbolic, at least two eigenvalues,
say $\theta_1$ and $\theta_2$, are distinct. Hence there exist
$v_1, v_2, \dots, v_l \in V_{\lambda_j}^\mathbb C$ (the
complexification of $V_{\lambda_j}$)
 such that $\overline\tau(v_i) = \theta_i v_i + x_i$, where $x_i
 \in
 \N_{j-1}^\mathbb C$, for all $1 \leq i \leq l$. We note that
 $v_1$ and $v_2$ are linearly independent since $\theta_1$ and $\theta_2$ are
 distinct.
We write $v_i = \sum_{\alpha \in S_{\lambda_j}}\! a_\alpha^i
\alpha$ where $a_\alpha^i \in \mathbb C$ for all $\alpha \in
S_{\lambda_j}$ and $1 \leq i \leq l$.
 As $v_1$ and $v_2$ are linearly independent and $\alpha \beta \in E$
for all $\alpha, \beta \in S_{\lambda_j}$, we have $[v_1, v_2]
\neq 0$ in $\N^\mathbb C$. Hence $[v_l,[ \cdots,[v_2, v_1]\cdots]]
 \neq 0$ in $\h_k$.
 Let $x =  [v_l,[ \cdots,[v_2, v_1]\cdots]]$. By considering the
complexification of $\tau$ and $\overline \tau$ we have
$\tau^\mathbb C(x) = (\prod_{i=1}^{l}\theta_i) x  + y$, where $y$
belongs to the complexification of $[\N_j,[\N_j,\cdots [\N_j,
 \N_{j-1}]\cdots]]_{(l -1) \mbox{ times }}$ which we denote  by
$W'$. We note that $\prod_{i=1}^{l}\theta_i = \pm1$ and $x \notin
W'$.
 Hence we
have an eigenvalue $\pm 1$ for the induced action of $\tau$ on
$(\N_j^l)^\mathbb C /
W'$ which is  a contradiction, since by assumption $\tau$ is
hyperbolic.
 This shows
 that $\alpha \beta$ is not an edge for all $\alpha, \beta \in
 S_{\lambda}$,
 where
 $|S_\lambda| = l, 1 \leq l \leq k$. This completes the proof of the theorem.
\hfill$\square$

\begin{examples}
{\rm i) Let $(S, E)$ be a complete graph i.e. $\alpha \beta \in E$
for all $\alpha, \beta \in S$. Then the corresponding $k$-step
nilmanifold admits an Anosov automorphism if and only if $|S| > k$.

ii) Let $(S , E)$ be a cycle on $4$ vertices. The corresponding
$k$-step nilmanifold admits an Anosov automorphism for all $k \geq
2$.
 In particular, we
get an example of $20$-dimensional
  $3$-step nilmanifold admitting an Anosov automorphism.

iii) A  complete bipartite graph $(S , E)$ is a graph where
$S$ is a disjoint
union of two subsets $S_1$ and $S_2$, each containing at least two
elements,  and $E = \{ \alpha \beta : \alpha
\in S_1, \beta \in S_2 \}$. In this case
$S_1$ and $S_2$ are the coherent components.
 Hence the $k$-step nilmanifold associated with a
complete bipartite graph admits an Anosov automorphism for all $k
\geq 2$. In particular, if we choose $S_1$ and $S_2$ such that $|S_1| =
 m$ and $|S_2| = n$ we get an example of $l$-dimensional $3$-step
nilmanifold admitting an Anosov automorphism, where $l =
 m (n-1)^2 - \frac{(n-2) (n-1)  m}
{2} + n  (m-1)^2 - \frac{(m-2) (m-1) n}  {2} + 2mn$.

iv) Let $(S , E)$ be a ``magnet'' graph with core $C$ i.e. $C$ is a
subset of $S$ such that its complement in $S$ contains at least two
elements and $E = \{ \alpha \beta : \alpha \in C, \beta \in S , \alpha
\neq \beta \}$. The $k$-step nilmanifold associated with $(S , E)$
admits an  Anosov automorphism if and only if  $k < |C|$.

}

\end{examples}
\section{Nonexistence of Anosov automorphisms on certain 2-step
nilmanifolds}

In this section we  prove  some results on nonexistence of
 Anosov automorphisms on
certain nilmanifolds. Let $\N_\mathbb Q$  be the
 $2$-step nilpotent Lie algebra over $\q$,
 associated
to the graph $(S, E)$.
 Let $X = [\alpha, \beta] + [\gamma, \delta]$, where $\alpha, \beta,
\gamma, \delta$ are distinct vertices in $S$ such that $\alpha\beta, \gamma\delta ,
\alpha \gamma , \alpha \delta \in E$.
 Let $\h_\mathbb Q$ denote the quotient
$\N_\mathbb Q/\big{<}X\big{>}$ where $\big{<}X\big{>}$ is the one-dimensional
 subspace spanned by $X$. Let $\h = \N / \big{<}X\big{>}$.
 It was proved in  \cite{De-Des} that if the graph $(S,E)$ is a complete
graph
 (i.e. $\alpha \beta \in E$ for all  $\alpha, \beta  \in S$), then
 $\h_\mathbb Q$ does not admit  a hyperbolic
automorphism  whose characteristic polynomial has integer coefficients and unit
constant term (see Theorem 3.2 of \cite {De-Des}).
 We prove a  similar result for an arbitrary graph.

\begin{theorem}\label{q2step}
The $2$-step nilpotent Lie algebra $\h_\mathbb Q$, defined as above,
 does not admit a hyperbolic automorphism
whose characteristic polynomial has integer coefficients and unit constant term.
\end{theorem}

\begin{notation}\label{def-G}
{\rm We recall that $\N = V \oplus W$ (see  \S\ref{k-stepdef}).  We decompose $\h$
as $\h = V \oplus W'$, where $W' = W/ \big{<} X \big{>}$.
%Let $\overline T$ be the
%subgroup of $\Aut(\h)$ defined as $\overline T = \{ \tau \in \Aut(\h) :\tau(V) =
%V\}$.
 Let $\overline G$ be the  subgroup of $GL(V)$ consisting of all restrictions,
$\tau|V$, such
 that $\tau \in \Aut(\h)$ and $\tau(V) = V$. Let $G$ be  the subgroup of $GL(V)$
 consisting of all
 restrictions, $\tau|V$, such that $\tau \in \Aut(\N)$
  and $\tau(V) = V$.
  It can be seen that  subgroups $G$ and  $\overline G$ of $GL(V)$ are  Lie
subgroups. Let $\G$ (resp. $\overline \G$) be the  Lie algebra of $G$ (resp.
 $\overline G$). Let $G^0$ (resp. $\ol G^0$) be the connected
 component of the identity in $G$ (resp. in $\ol G$). Let $\D$ (resp. $\overline \D$)
  be the Lie subalgebra of $\G$ (resp. of $\overline \G$) consisting of all
 endomorphisms in $\G$ (resp. in $\ol \G$) that are represented by
 diagonal matrices with respect
 to the basis $S$. Note that $\ol \D$ consists of all the
 endomorphisms in $\D$  which are contained in $\ol \G$.
 For $\eta, \zeta \in S$,
 let $E_{\eta \zeta}$ be the
 element of $\End (V)$ such that $E_{\eta \zeta}(\zeta)=\eta$ and
 $E_{\eta \zeta}(\xi)=0$ for all $\xi \in S$, $\xi \neq \zeta$.
}
\end{notation}

\begin{notation}
 {\rm Recall that $X = [\alpha, \beta] + [\gamma, \delta]$,
 where $\alpha, \beta,
 \gamma, \delta$ are distinct vertices in $S$ such that
  $\alpha\beta, \gamma\delta,
\alpha\gamma, \alpha \delta \in E$. Let $S' = \{ \alpha, \beta,
\gamma, \delta \}$.
 Let $W_{\phi \psi}^{\phi' \psi'}$  denote the subspace of $\End(V)$
spanned by $E_{\phi \psi}$ and $E_{\phi' \psi'}$ where
 $\{\phi, \psi, \phi',  \psi' \} = S'.$
}
\end{notation}

\begin{proposition}\label{span}
The Lie algebra  $\ol \G$, defined as above, is spanned by
 $\ol \D$,   $W_{\alpha  \delta}^{\gamma \beta} \cap
\ol \G$,  $W_{\beta \gamma}^{\delta \alpha} \cap \ol \G$,
$W_{\alpha \gamma}^{\delta \beta} \cap \ol \G$, $W_{\gamma
\alpha}^{\beta \delta} \cap \ol \G$, and the elements of $\ol \G$
of the following type: (i) $E_{\eta \zeta}$, where $\eta \neq
\zeta$, $\eta, \zeta \notin S'$, (ii) $E_{ \eta \zeta}$, where
$\eta \in S'$ and $\zeta \notin S'$, (iii) $E_{ \eta \zeta}$,
where $\eta \notin S'$ and $\zeta \in S'$, (iv) $E_{ \alpha
\beta}$, $E_{\gamma \delta}$, $E_{\beta \alpha}$, $E_{\delta
\gamma}$.

\end{proposition}

\proof Let $Y \in \ol \G$. Then it can be expressed as $Y = Y_0 +
\sum_{\eta, \zeta \in S, \eta \neq \zeta} a_{\eta \zeta} E_{\eta
\zeta}$, where $Y_0 \in \D$ (see Notation \ref{def-G}) and
$a_{\eta \zeta}\in \R$.  By using the fact that $E_{\zeta \zeta}
\in \overline \G$ for all $\zeta \notin S'$,
 we  observe that  $a_{\eta \zeta}
E_{\eta \zeta}$ is contained in $\ol \G$ for all $\eta, \zeta \notin S'$ (see the
proof of Proposition $3.1$ in \cite{D-M}). We note that $E_{\alpha \alpha} +
E_{\gamma \gamma}$ , $E_{\beta \beta} + E_{\delta \delta}$ , $E_{\beta \beta}+
 E_{\gamma \gamma}$ are contained in
$\ol \G$. Since $[E_{\zeta \zeta}, [E_{\alpha\alpha} + E_{\gamma
\gamma}, Y]]$ and $[E_{\zeta \zeta}, [E_{\zeta
\zeta},[E_{\alpha\alpha} + E_{\gamma \gamma}, Y]]]$ are in $\ol
\G$ for $\zeta \notin S'$, $a_{\zeta \alpha}E_{\zeta \alpha} +
a_{\zeta \gamma} E_{\zeta \gamma}$ and $a_{\alpha \zeta} E_{\alpha
\zeta} + a_{\gamma \zeta}E_{\gamma \zeta}$ are contained in $\ol
\G$. Now as $E_{\alpha \alpha} + E_{\delta \delta} \in \ol \G$, we
have $a_{\zeta \alpha} E_{\zeta \alpha}$ and $a_{\alpha
\zeta}E_{\alpha \zeta}$ are in $\ol \G$. Similarly we can see that
$a_{\eta \zeta}E_{\eta \zeta} \in \ol \G$, for all $\eta \in S'$
and $\zeta \notin S'$; and $a_{\eta \zeta} E_{\eta \zeta} \in \ol
\G$, for all $\eta \notin S'$ and $\zeta \in S'$. We also have
$[E_{\beta \beta}+ E_{\gamma \gamma} ,[E_{\beta \beta} + E_{\delta
\delta},[E_{\alpha \alpha} + E_{\gamma \gamma}, Y]]] \in \ol \G$.
This shows that $Z = a_{\gamma \delta}E_{\gamma \delta} + a_{\beta
\alpha}E_{\beta \alpha} -a_{\delta \gamma}E_{\delta \gamma} -
a_{\alpha \beta}E_{\alpha \beta} \in \ol \G$. Also $[ E_{\beta
\beta}+ E_{\gamma \gamma}, Z] \in
 \ol \G$. Therefore we get $a_{\gamma \delta}E_{\gamma \delta}
+ a_{\beta \alpha}E_{\beta \alpha} +a_{\delta \gamma}E_{\delta \gamma} +a_{\alpha
\beta}E_{\alpha \beta} \in \ol \G$. Hence $a_{\gamma \delta}E_{\gamma \delta} +
a_{\beta \alpha}E_{\beta \alpha}$ and $a_{\delta \gamma}E_{\delta \gamma} +a_{\alpha
\beta}E_{\alpha \beta}$
 are contained in $\ol \G$.
Since $[E_{\alpha \alpha} + E_{\gamma \gamma}, a_{\gamma \delta} E_{\gamma \delta} +
a_{\beta \alpha}E_{\beta \alpha}] \in \ol \G$, we have $a_{\gamma \delta}E_{\gamma
\delta} - a_{\beta \alpha}E_{\beta \alpha} \in \ol \G$ and hence $a_{\gamma
\delta}E_{\gamma \delta}\in \ol \G$. Hence we have proved that if $a_{\gamma \delta}
\neq 0$ then  $E_{\gamma \delta}\in \ol \G$. Similarly it can be proved that
$E_{\beta \alpha} \in \ol \G$ if $a_{\beta \alpha} \neq 0$, $E_{\alpha \beta} \in
\ol \G$ if $a_{\alpha \beta} \neq 0$, and $E_{\delta \gamma} \in \ol \G$ if
$a_{\delta \gamma} \neq 0$. Now $Z' = [E_{\alpha \alpha} + E_{\gamma \gamma},
[E_{\beta \beta} + E_{\delta \delta},[E_{\alpha \alpha} + E_{\gamma \gamma}, Y]]]
\in \ol \G$, and hence
  $[E_{\beta \beta}
+ E_{\delta \delta},[E_{\alpha \alpha} + E_{\gamma \gamma}, Y]] + Z' \in \ol \G$. As
$a_{\alpha \beta} E_{\alpha \beta}$ and $a_{\gamma \delta}E_{\gamma \delta}$ are
contained in $\ol \G$, we have $a_{\alpha \delta}E_{\alpha \delta}
 + a_{\gamma \beta}E_{\gamma \beta} \in \ol \G$. Similarly we can
prove that $a_{\beta \gamma} E_{\beta \gamma} + a_{\delta \alpha}
 E_{\delta \alpha} \in \ol \G$.
As $Y \in \ol \G$, by above observations we have $Z'' = Y_0 + a_{\alpha
\gamma}E_{\alpha \gamma} + a_{\gamma \alpha}E_{\gamma \alpha} + a_{\beta
\delta}E_{\beta \delta} + a_{\delta \beta}E_{\delta \beta} \in \ol \G$. Considering
 the element $[E_{\alpha \alpha} + E_{\delta \delta}, Z'']$
we prove that $a_{\alpha \gamma}E_{\alpha \gamma} + a_{\delta \beta}E_{\delta \beta}
\in \ol \G$ and $a_{\gamma \alpha}E_{\gamma \alpha} + a_{\beta
\delta}E_{\beta\delta}$ are in $\ol \G$. Hence we have now $Y_0 \in \ol \G$. As
$Y_0$ is in $\D$, $Y_0 \in \ol \D$. Hence we have proved our claim  that $\ol \G$ is
spanned by $\ol \D$,  $W_{\alpha  \delta}^{\gamma \beta} \cap \ol \G$,
 $W_{\beta \gamma}^{\delta \alpha} \cap
\ol \G$, $W_{\alpha \gamma}^{\delta \beta} \cap \ol \G$, $W_{\gamma \alpha}^{\beta
\delta} \cap \ol \G$,
 and the elements of $\ol \G$ of the type (i)-(iv) as stated.
\hfill$\square$

\begin{proposition}\label{lemma1}
Any automorphism $\ol T$ in  $\ol G^0$  is induced by  an automorphism $T$ in $G^0$,
with  $T\big{(}\big{<} X \big{>}\big{)}  = \big{<}X \big{>}$.
\end{proposition}

\proof We will prove the following: If the element from  the type (i)-(iv)
 in the statement of Proposition \ref{span}, considered as an
element of $\End(V)$, is in the Lie algebra $\ol \G$; then that element is in $\G$
(see Notation \ref{def-G}). We will also prove that  $W_{\alpha  \delta}^{\gamma
\beta} \cap \ol \G$,
   $W_{\beta \gamma}^{\delta \alpha} \cap
\ol \G$, $W_{\alpha \gamma}^{\delta \beta} \cap \ol \G$, and $W_{\gamma
\alpha}^{\beta \delta} \cap \ol \G$ are contained in $\G$. Let $I$ denote the
identity transformation in $GL(V)$.

(i) If $E_{\eta \zeta} \in \ol \G$, where $\eta \neq \zeta$, $\eta, \zeta \notin
S'$; then $\tau = I + E_{\eta \zeta} \in \ol G$. We will prove that $\Omega'(\eta)
\subset \Omega(\zeta)$ (see \S\ref{auto}). Suppose that $\xi  \notin \Omega(\zeta)$.
Then we have $[\xi, \zeta] = 0$. As $\tau \in \ol G$, it is the restriction of a Lie
automorphisms of $\h$, and hence we get that $[\tau(\xi), \tau(\zeta)] = 0$ in $\h$.
Therefore we have $[\xi, \zeta + \eta] = c ([\alpha, \beta] + [\gamma, \delta])$,
where $c \in \R$. As $[\xi, \zeta] = 0$, we have $[\xi, \eta ] = c ([\alpha, \beta]
+ [\gamma, \delta])$ and hence $c = 0$.  Hence $\xi \notin \Omega'(\eta)$. This
shows that
 $\Omega'(\eta) \subset \Omega(\zeta)$.
  Hence  $E_{\eta \zeta} \in \G$,
 where $\eta \neq \zeta$, $\eta, \zeta \notin S'$ (see Proposition $4.1$ in \cite{D-M}).

(ii) Consider the element $E_{\eta \alpha}$, where $\eta \notin S'$. Suppose
$E_{\eta \alpha} \in \ol \G$ Let $\tau = I + E_{\eta \alpha}$.
 Let $ \zeta \notin \Omega(\alpha)$. Then $[\zeta, \alpha] = 0$. As
$\tau \in \ol \G$, we have $[\zeta, \alpha + \eta] = 0$ in $\h$. By the similar
argument as above we have $[\zeta , \eta] = 0$. Hence we have $\Omega'(\eta) \subset
\Omega(\alpha)$. Hence $E_{\eta \alpha} \in \G$.

By  similar arguments we can prove our claim for the elements of the type (ii),
(iii) and (iv).

We will prove that any element of  $W_{\alpha \delta}^{\gamma
\beta} \cap
  \ol \G$, considered as an element of $\End(V)$, is contained in $\G$.
   Suppose now the linear combination of $E_{\alpha
\delta}$ and $E_{\gamma \beta}$, say $a E_{\alpha \delta} +
 b E_{\gamma \beta}$, is in $\ol \G$. Then $\tau = I + t(E_{\alpha \delta} +
 b E_{\gamma \beta})\in \ol G, t \neq 0$. We show that the subspace  of
$\wedge^2 V$ spanned by the set of all $\zeta \wedge \eta$ such that $\zeta \eta$ is
not an edge, say $W'$, is $\wedge^2 \tau$-invariant. Let $\zeta, \eta \in S$ such
that $\zeta \neq \eta$ and $\zeta \eta$ is not an edge. If neither of $\zeta$ and
$\eta$  is contained in $\{\delta, \beta\}$, then $\wedge^2 \tau(\zeta \wedge \eta)
=  \zeta \wedge \eta$. If $\zeta = \delta$ and $\eta \neq \beta$ then as $\tau \in
\ol G$, we have $\tau[\zeta, \eta] = 0$ in $\h$, and hence $[ ta\alpha + \delta,
\eta] = c ([\alpha, \beta] + [\gamma, \delta])$ in $\N$, where $c \in \R$. Therefore
either $a = 0$ or $\alpha \eta$ is not an edge. In both the cases we have $\wedge^2
\tau(\zeta \wedge \eta) \in W'$. Similarly if $\zeta = \beta$ and $\eta \neq \delta$
we are through. Now if $\zeta = \delta$ and $\eta = \beta$, we have $\tau[\delta,
\beta] = 0$ in $\h$. Hence $[ta \alpha + \delta, tb\gamma + \beta] = c ([\alpha,
\beta] + [\gamma, \delta])$ in $\N$, where $c \in \R$.
 As $\delta \beta$ is not an
edge and $\alpha \gamma$ is an edge,
 we have $t^2ab[\alpha, \gamma] + ta [\alpha ,
\beta] + tb [\delta , \gamma] =
 c ([\alpha, \beta] + [\gamma, \delta])$ in $\N$, and hence $ab = 0$
and $a = -b$. Therefore $\tau = I$. Thus we have proved that $W'$ is $\wedge^2
\tau$-invariant. Therefore $\tau \in G$, and hence $a E_{\alpha \delta} +
 b E_{\gamma \beta} \in \G$.

Similarly we see that our claim holds for the elements of $W_{\beta \gamma}^{\delta
\alpha} \cap \ol \G$, $W_{\alpha \gamma}^{\delta \beta} \cap \ol \G$, and $W_{\gamma
\alpha}^{\beta \delta} \cap \ol \G$.

By using the above argument, Proposition \ref{span} and
 Theorem 2.10.1 in \cite{V}, we
see that there exists an open neighbourhood $U$ of $I$ in $\ol G^0$
 such that any automorphism contained in $U$ can be lifted to an
automorphism of $\N$ which keeps an ideal $\big{<} X \big{>}$ invariant. Hence any
automorphism $\ol T$ in $\ol G^0$ can be lifted to an automorphism $T$ in $G^0$ such
that $T(\big{<} X \big{>}) = \big{<} X \big{>}$ (use Proposition 3.18 in \cite{W}).
\hfill$\square$

\bigskip

\noindent {\bf Proof of  Theorem \ref{q2step}: }
 Suppose $\ol \theta \in \Aut(\h_\mathbb Q)$ is a hyperbolic automorphism
such that its  characteristic polynomial has integer coefficients and unit constant
term. Since $\Aut(\h)$ has  finitely many connected components,
 by replacing $\ol \theta$ by its suitable power we may assume that
$\ol \theta$ is contained in the connected component of the identity in $\Aut(\h)$.
By Proposition
 \ref{lemma1}, we see that there exists an automorphism
$\theta$ contained in the connected component of the identity of $\Aut(\N)$ such
that its  characteristic polynomial has integer coefficients and unit constant,
$\theta(X) = X$, and $\theta$ has an eigenvalue $1$ of multiplicity $1$. We can
assume that the matrix of $\theta$ with respect to the basis $S \cup E$ is an
integer matrix.

We have $\theta(\N_j) = \N_j$ for each  $j = 1,\dots,k$ where
$\N_j = (\oplus_{i\leq j}\,V_{\lambda_i})\oplus W$ and $\lambda_1,
\dots ,\lambda_k$ is  an enumeration of $\Lambda$ such that for
all $j= 1,\dots, m$, $\bigoplus_{i\leq j}\,V_{\lambda_i}$ is
invariant under the  action of $G^0$ (see \S \ref{2-step}). Let
$\pi_j : \N_j \rightarrow V_{\lambda_j}$ denote the canonical
projection for each $j = 1,\dots,k$. Let $\theta_{\lambda_j} :
V_{\lambda_j} \rightarrow V_{\lambda_j}$ be given by
$\theta_{\lambda_j} = \pi_j \circ \theta$.

We have $W = \sum_{\lambda, \mu \in \Lambda} [V_{\lambda}, V_{\mu}]$.
 All the eigenvalues of $\theta$ on $W$  are pairwise products
  of the eigenvalues
  on $V_{\lambda}$'s. Also $a_\lambda a_\mu$
 occurs as an eigenvalue of $\theta|W$ if and only if there exists
$\zeta \in S_\lambda$ and $\eta \in S_\mu$ (notation are as before) such that $\zeta
\eta$ is an  edge, where $a_\lambda$ is an eigenvalue of $\theta_{\lambda}$ and $
a_\mu$ is an eigenvalue of $\theta_{\mu}$ As $1$ is an eigenvalue of $\theta$, there
exist $\lambda$ and $\mu$ in $\Lambda$ such that $\zeta \eta$ is an  edge for $\zeta
\in S_{\lambda}$ and $\eta \in S_{\mu}$, and $a_\lambda a_\mu = 1$,
 $a_\lambda$ and  $a_\mu$
being eigenvalues of $\theta_{\lambda}$ and $\theta_{\mu}$ respectively.

We will prove that $\lambda = \mu$. Suppose that $\lambda \neq \mu$.
 Let  $a_{\lambda}'$ be a conjugate of $a_{\lambda}$ over $\mathbb Q$,
 $a_{\lambda}' \neq a_{\lambda}$.  Then $a_{\mu} =
a_{\lambda}^{-1}$ and $a_{\lambda'}^{-1}$ are conjugates over $\mathbb Q$. The
minimal polynomial of $a_{\mu}$ over $\mathbb Q$ divides the characteristic
polynomial of $\theta_{\mu}$. Hence $a_{\lambda'}^{-1}$ occurs as an eigenvalue of
$\theta_{\mu}$. As by our assumption $\zeta \eta$ is an edge for all $\zeta \in
S_{\lambda}$ and $\eta \in S_{\mu}$, $a_{\lambda'}a_{\lambda'}^{-1}$ occurs as an
eigenvalue of $\theta$,  where $a_{\lambda} \neq a_{\lambda}'$, and hence we arrive
at a contradiction as the multiplicity of the eigenvalue 1 is 1. Therefore $\lambda
= \mu$. Hence there exists $\lambda \in \Lambda$ such that the restriction of a
graph $(S, E)$ on $S_{\lambda}$ is complete  and $a_{\lambda} a_{\lambda}' =1$ for
the some eigenvalues $a_{\lambda}$ and $a_{\lambda}'$ of $\theta_{\lambda}$ and $X
\in [V_\lambda, V_\lambda]$, which is not possible (by Theorem 3.2 of \cite{D}).
This completes the proof of the theorem. \hfill$\square$

%\midskip

\begin{remark}
{\rm Let $H$ be the simply connected nilpotent Lie group
corresponding to the Lie algebra $\h$. Let $\Gamma$ be a lattice
in $H$ corresponding to $\h_\mathbb Q$ (see \cite{D}). Then
Theorem \ref{q2step} shows that the nilmanifold $\Gamma \backslash
 H$ does not admit an Anosov automorphism. }
\end{remark}

\section{Nonexistence of Anosov automorphisms on $3$-step
nilmanifolds}

In this section  we  study  quotients of certain $3$-step
nilpotent
 Lie algebras. Let
$(S,E)$ be a graph. Let $V$ be a vector space over $\mathbb Q$
with a basis as $S$.
 Let $\F_3(V)$ denote the free 3-step nilpotent Lie algebra over
  the rationals on $V$. Let $\Q = \F_3(V)/ \I $, where $\I$ is an ideal
of $\F_3(V)$ generated by the elements $ [\alpha , \beta]$,
 where $\alpha, \beta \in
S$ and $\alpha \beta$ is not an edge.  We decompose
 $\Q$ as $\Q = V \oplus W \oplus V_3$, where $V_3$ is the space
spanned by all $[\alpha, [\beta, \gamma]]$ such that $\alpha, \beta,
\gamma\in S$,
and $\beta \gamma$ is an edge. Let $X$ be a nonzero vector in $V_3$. Let
  $ \M  = \Q /\big{<}X\big{>}$, where $\big{<}X\big{>}$ denotes
  the ideal  generated by $X$ in $\Q$ which is a one-dimensional
subspace spanned by $X$ in $\Q$.

\begin{proposition}
Let $\overline{\theta}$  be an automorphism of $\M$. Then there
exists $\theta$, an automorphism of $\Q$ such that $\theta (X) = c
X$, where c is a nonzero rational. Furthermore $\theta$ induces an
automorphism $\theta'$ of $\M$ such that both $\overline{\theta}$
and $\theta'$ have the same eigenvalues and also we have
$\theta(V) = V, \theta(W) = W$ and $\theta(V_3) = V_3$.
\end{proposition}
\proof Consider the linear endomorphism  $T$ of $ V$ defined by
$T(\alpha) = \pi
(\overline{\theta}(\overline{\alpha}))$ for all $\alpha \in S$,
where $\pi$ is a
natural projection of $\M$ onto $V$ with respect to the decomposition
of $\M$ as $\M
= V \oplus W \oplus V_3/ \big{<}X\big{>}$,
  and
$\overline{\alpha}$ denotes the coset in $\M$ represented by $\alpha$. Note that $T$
is an automorphism of $V$.
 For, if $T(v) = 0, v \in V$, then $\overline{\theta}(v) \in [\M,
\M].$ But as $\overline{\theta}$ is an automorphism of $\M$,
we must have $v=0$.

We then have a Lie algebra  automorphism $\overline{T}$ of
 $\F_3(V)$ such that $\overline{T}|_V = T$.
We prove that $\overline{T}(\I) = \I$. Suppose $\alpha, \beta$ are in $ S$ such that
$ \alpha \beta $ is not an edge. By definition, $\overline{T}[\alpha, \beta] =
[\pi(\overline{\theta} (\overline {\alpha})),  \pi(\overline{\theta} (\overline
{\beta}))]$. Now as $\alpha \beta$ is not an edge, we have $[\overline {\alpha},
\overline {\beta}] = 0$. Hence $[ \overline {\theta}(\overline{\alpha}) , \overline
{\theta}(\overline{\beta})] = 0$ in $\M$. It implies that $[\pi(\overline{\theta}
(\overline {\alpha})), \pi(\overline{\theta} (\overline {\beta}))] \in V_3.$  Hence
we have $[\pi(\overline{\theta} (\overline {\alpha})), \pi(\overline{\theta}
(\overline {\beta}))] = 0$ in $\Q$. That means in $\F_3(V)$, we have
$[\pi(\overline{\theta} (\overline {\alpha})), \pi(\overline{\theta} (\overline
{\beta}))] \in \I$. Thus we have proved $\overline{T}(\I) = \I$ and hence
  $\overline{T}$ factors through an automorphism of $\Q$.
Let $\theta$ denote the automorphism of $\Q$ induced by
$\overline{T}$.

We claim that $\theta(X) \in \big{<}X \big{>}.$ We note that
$\overline{\theta}(\overline {Y}) = \overline {\theta(Y)}$ in $\M$, for all $Y \in
V_3$, where bar is taken to denote the elements in $\M$ represented by elements in
$\Q$. Now $\overline{\theta}(\overline{X}) = \overline{0}$ in $\M$, as
 $\overline{\theta}$ is an automorphism of $\M$. This implies that
$\theta(X)\in\big{<}X \big{>}$ in $\Q$. Hence we have $\theta(X) =
c X$, where $c$ is a nonzero rational. By definition of $\theta$,
$\theta(V) = V, \theta (W) = W$, and $\theta (V_3) = V_3$. Let
$\theta'$ be the  automorphism of $\M$ induced by $\theta$. Then
 both $\overline{\theta}$ and $\theta'$ induce the same linear
endomorphism on $\M / [\M, \M]$. Hence  $\overline{\theta}$ and $\theta'$ have same
eigenvalues. \hfill$\square$

\begin{notation}\label{E}
{\rm Let $\{S_\lambda\}_{\lambda \in \Lambda}$ denote the family of
coherent components (see \S \ref{2-step}) of the graph $(S, E)$.
 Let $\cal E$ be the set of all unordered pairs $\lambda \mu$
 with  $\lambda, \mu \in \Lambda$, such that
 $\alpha \beta \in E$ for $\alpha \in S_\lambda$ and $\beta \in
S_\mu$.
 We recall that $V_\lambda$
denotes the subspace of $V$ spanned by $S_\lambda$,
 $\lambda \in \Lambda$.
 Let $\lambda_1, \dots ,\lambda_k$ be an enumeration of $\Lambda$ such that
assertion~(ii) of Theorem~\ref{thm} holds.
Let  $\N_j = (\oplus_{i\leq j}\,V_{\lambda_i})\oplus W$ for $j=1,\dots,k$.}
\end{notation}

\begin{theorem}\label{last}

The Lie algebra $\M$ does not admit a hyperbolic automorphism whose
 characteristic
polynomial has integer coefficients and unit constant term.
\end{theorem}
\proof Let $\overline{\theta}$ be a hyperbolic
 automorphism of $\M$ such that the
characteristic polynomial of $\overline{\theta}$
   has integer coefficients and unit constant term. Let $\theta$ be an
automorphism of $\Q$ as obtained in the previous proposition. Let $\tau$ be an
automorphism of $\N$ induced by $\theta$.
 As $\overline{\theta}$ is a hyperbolic
 automorphism, and $\theta(X) = X$, 1 is an eigenvalue of $\theta$ of
 multiplicity 1.
Since the characteristic polynomial of $\overline{\theta}$ has integer
coefficients
and unit constant term, we may assume that the matrix of $\theta$ has all integer
entries (by replacing $\theta$ by some power of $\theta$ if necessary.) As
$\mbox{Aut}(\N)$ has finitely many components,  we may assume that
$\theta$
 lies in
the connected component of the identity of $\mbox{Aut}(\Q)$ and $\tau$
lies  in the
connected component of the identity of $\mbox{Aut}(\N)$. Hence
 $\theta(\N_j) =
\N_j$ (see Lemma \ref{lem:det}).

Let $\pi_{\lambda_j} : \N_j \rightarrow V_{\lambda_j}$ denote the canonical
projection. Let $\theta_{\lambda_j}$ be an endomorphism of $
V_{\lambda_j}$
 defined
 by $\theta_{\lambda_j}(v) = \pi_{\lambda_j}(\theta(v))$, for all $v \in
 V_{\lambda_j}$. We note that $\theta_{\lambda_j}$ is an automorphism of
 $V_{\lambda_j}$.

All the eigenvalues of $\theta$ restricted on $V_3$
 are of the following type:

i) $\delta_\lambda \delta_\mu \delta_\nu$, where $\delta _\lambda,
 \delta_\mu,$ and $ \delta_\nu$ are the eigenvalues of $\theta_\lambda,
 \theta_\mu$, and $\theta_\nu$ respectively and $\mu \nu \in \cal E$
(see Notation \ref{E}).

ii) $\delta_\lambda \delta_\lambda' \delta_\lambda''$, where $\delta_\lambda,
\delta_\lambda', \delta_\lambda''$ are the eigenvalues of $\theta_\lambda$ and the
restriction of $(S, E)$ to $S_\lambda$ is a complete graph.

As $\theta|_{V_3}$ has an eigenvalue 1,  we have
 the following two cases:

\noindent {\bf Case(i):} Suppose  $\delta_\lambda \delta_\mu
\delta_\nu = 1$,
 where
$\delta _\lambda, \delta_\mu,$ and $ \delta_\nu$ are the eigenvalues of
$\theta_\lambda,\theta_\mu$, and $\theta_\nu$ respectively and $\mu \nu \in \cal E$.
Now as  $\mu \nu \in \cal E$, $\delta_\mu \delta_\nu$ occurs as an
 eigenvalue of
$\theta|_W$. Thus we have an invertible matrix, say $A$,
 with integer entries such
that
 $\delta_\mu \delta_\nu$ is an eigenvalue of $A$. Also we have an
invertible matrix, say $B$, with integer entries having
 $\delta_\lambda = (\delta_\mu \delta_\nu)^{-1}$ as an
 eigenvalue. Hence there exists an eigenvalue of the type $\delta_\mu'
 \delta_\nu'$
of $A$ and $\delta_\lambda' = (\delta_\mu' \delta_\nu')^{-1}$ of $B$ such that
$\delta_\mu', \delta_\nu'$ and $\delta_\lambda'$ are the
 eigenvalues of $\theta_\mu,
\theta_\nu$ and $\theta_\lambda$ respectively,
 and $\delta_\mu' \delta_\nu' \neq
\delta_\mu \delta_\nu$.  This contradicts the fact that
 the multiplicity of the eigenvalue 1 is 1.

\noindent {\bf Case(ii):}
  Suppose
 $\delta_\lambda \delta_\lambda' \delta_\lambda'' = 1$
  for some $\lambda \in \Lambda $ such that the restriction of $(S,E)$
on $S_\lambda$ is complete and $\delta_\lambda, \delta_\lambda',
 \delta_\lambda''$
are the eigenvalues of $\theta_\lambda$. If $\delta_\lambda =
\delta_\lambda' = \delta_\lambda''$ then $\delta_\lambda^3 = 1$,
 which is a contradiction since $\overline{\theta}$ is hyperbolic.
 Hence we may assume that $ \delta_\lambda'\neq
\delta_\lambda''$. Let $V_\lambda^{\mathbb{C}}$ denote the
complexification of $V_\lambda$.
 Suppose that
$Y, Y', Y'' \in V_\lambda^{\mathbb{C}}$ are eigenvectors corresponding to the
eigenvalues $\delta_\lambda, \delta_\lambda', \delta_\lambda''$. We
consider
the
complexification of $\Q$ and $\theta$ also. As $\delta_\lambda' \neq
\delta_\lambda'',$ $Y'$ and $Y''$ are linearly independent.
Hence $[Y, [Y' , Y'']]
\neq 0$ in $\Q$ and also we have $\theta[ [Y, [Y' , Y'']] =
[Y, [Y' , Y'']]$.  If all
the $\delta_\lambda, \delta_\lambda',
 \delta_\lambda''$ are distinct, then consider $[Y', [Y , Y'']]$, which is
 an eigenvector corresponding to the eigenvalue 1. Also
$[Y', [Y , Y'']]$ and $[Y, [Y' , Y'']]$ are linearly independent. This is not
possible as the eigenvalue 1 has multiplicity 1. If $\delta_\lambda =
\delta_\lambda''$ and $Z \in W^{\mathbb{C}}$
 is an eigenvector of $\theta$ corresponding to the eigenvalue
$\delta_\lambda^{2}$, then consider $[Y', Z]$, which is an eigenvector
 of $\theta$ corresponding to the eigenvalue 1. Also $[Y', Z]$ and
$[Y,[Y', Y'']$ are linearly independent. This is a contradiction. Similarly we get a
contradiction if $\delta_\lambda = \delta_\lambda'$. This proves the theorem.
\hfill$\square$

\begin{remark}
{\rm Theorem \ref{last} shows that
 the nilmanifold $\Gamma \backslash M$, where $M$ is the simply connected
 nilpotent Lie group corresponding to the Lie algebra $\M \otimes \R$ and
$\Gamma$ is a lattice in $M$ corresponding to $\M$, does not admit
an Anosov automorphism.
 In particular, a nilmanifold  $\Gamma \backslash M$, where $\Gamma$ corresponds to
  the rational  Lie algebra
   given by a quotient of free $3$-step nilpotent Lie
algebra by a one-dimensional ideal, does not admit an Anosov
automorphism. }

\end{remark}

\noindent {\it Acknowledgements.} I am very
 grateful to Prof. S. G. Dani for his valuable
help. I express my gratitude to  Prof. J. Lauret for his helpful
comments and suggestions. I would like to thank TWAS, Trieste-Italy
  and CIEM,
National University of Cordoba, Argentina for their support.

\end{document}